\documentclass[12pt]{article}
\usepackage{amssymb,amsmath,amsfonts,amsthm}
\usepackage{stmaryrd}
\usepackage{color}
\usepackage[latin1]{inputenc}
\usepackage{enumitem}
\usepackage{textcomp}
\usepackage{hyperref}

\usepackage[small,bf,hang]{caption} 

\usepackage{epsfig}
\usepackage{url,xargs}
\usepackage[font=small,skip=4pt]{caption}

\def\BBox{\kern  -0.2cm\hbox{\vrule width 0.2cm height 0.2cm}}

\newtheorem{theorem}{Theorem}[section]
\newtheorem{lemma}[theorem]{Lemma}

\newtheorem{definition}[theorem]{Definition}
\newtheorem{corollary}[theorem]{Corollary}
\newtheorem{proposition}[theorem]{Proposition}
\newtheorem{remark}[theorem]{Remark}
\newtheorem{observation}[theorem]{Observation}

\newtheorem{problem}[theorem]{Problem}

\graphicspath{{figures/}}

\setlist{noitemsep}

\newcommand{\N}{\mathbb{N}}
\newcommand{\R}{\mathbb{R}}
\newcommand{\Z}{\mathbb{Z}}

\DeclareMathOperator{\Me}{Me}

\newcommand{\ca}{\mathcal}

\title{\textbf{A unified approach to construct snarks with circular flow number 5}}

\author{Jan Goedgebeur$^1$\thanks{Supported by a Postdoctoral Fellowship of the Research Foundation Flanders (FWO)}, Davide Mattiolo$^2$, Giuseppe Mazzuoccolo$^3$\\ 
\footnotesize $^1$ Deparment of Applied Mathematics, Computer Science and Statistics\\[-3mm]
\footnotesize Ghent University, Krijgslaan 281 - S9, 9000 Ghent, Belgium.\\
\footnotesize $^2$ Dipartimento di Scienze Fisiche, Informatiche e Matematiche,\\[-3mm]
\footnotesize Universit\`{a} di Modena e Reggio Emilia, Via Campi 213/b, 41126 Modena, Italy.\\
\footnotesize $^3$ Dipartimento di Informatica,\\[-3mm]
\footnotesize Universit\`{a} degli Studi di Verona, Strada le Grazie 15, 37134 Verona, Italy.\\
}
\begin{document}
\maketitle

\begin{abstract}
The well-known 5-flow Conjecture of Tutte, stated originally for integer flows, claims that every bridgeless graph has circular flow number at most 5. It is a classical result that the study of the 5-flow Conjecture can be reduced to cubic graphs, in particular to snarks. However, very few procedures to construct snarks with circular flow number 5 are known.

In the first part of this paper, we summarise some of these methods and we propose new ones based on variations of the known constructions. Afterwards, we prove that all such methods are nothing but particular instances of a more general construction that we introduce into detail.

In the second part, we consider many instances of this general method and we determine when our method permits to obtain a snark with circular flow number 5.
Finally, by a computer search, we determine all snarks having circular flow number 5 up to 36 vertices. It turns out that all such snarks of order at most 34 can be obtained by using our method, and that the same holds for 96 of the 98 snarks of order 36 with circular flow number 5.

\end{abstract}



\section{Introduction}\label{intro}

In  this  paper,  we  will focus our attention on the following type of flows in graphs: let $r\geq2$ be a real number, a \textit{circular nowhere-zero $r$-flow}, $r$-CNZF from now on, in a graph $G$ is a flow in some orientation of $G$ such that the flow value on each edge lies in the interval $[1,r-1]$ and such that the sum of the inner and outer flow in every vertex is zero. 
The circular flow number of a graph $G$, denoted by $\Phi(G)$ and first introduced in~\cite{GTZ}, is the infimum of the real numbers $r$ such that $G$ has a circular nowhere-zero $r$-flow. 
It was conjectured by Tutte that any bridgeless graph has a (circular) nowhere-zero $5$-flow~\cite{Tu} and it is well-known that the study of Tutte's Conjecture can be restricted to the class of cubic graphs, in particular of \textit{snarks}, i.e.\ cyclically $4$-edge-connected  cubic  graphs  of  girth  at  least  $5$,  with  no  $3$-edge-colouring or equivalently, with no circular nowhere-zero $4$-flow. 

Flow numbers received much attention in the last decades (see Zhang's monograph~\cite{Z} for an overview), and some characterizations of graphs having a given circular flow number are known. In addition to classical results of Tutte, we can mention, among others,~\cite{SS} and~\cite{S} for regular graphs.
It is proven in~\cite{PaZh} that for any rational value $q$ in the interval $(4,5]$, there exists a graph with circular flow number exactly $q$, and an analogous result is proven in~\cite{LuSk} -- even if we restrict our attention to the class of snarks.
The question of whether the circular flow number of a snark could be exactly equal to 5 has the Petersen graph as a well-known positive answer. However no other examples were known at the time. Mohar asked in 2003~\cite{Mo} if the Petersen graph is the only possible one. In 2006 M\'a\v cajov\'a and Raspaud~\cite{MaRa06} gave a negative answer to Mohar's question by constructing an infinite family of snarks with circular flow number 5. More recently, Esperet, Tarsi and the third author~\cite{EMT16}, extending the method proposed in~\cite{MaRa06}, constructed a larger class of snarks with circular flow number $5$ and, among other results, show that deciding whether a given snark has circular flow number less than 5 is an NP-complete problem. Finally, by using methods introduced in~\cite{EMT16} another family of snarks having circular flow number $5$ was presented in~\cite{AKLM}. 

The main aim of this paper is to propose a unified and compact description of all such methods and the new ones introduced here. A summary of our results is Theorem~\ref{theo:summary}, which is the main result of the present paper.

Furthermore, using a computer search we will determine all snarks with circular flow number 5 up to order $36$. It turns out that all such snarks of order at most $34$ fit our description (and hence, more specifically, the description proposed in~\cite{EMT16}), and that the same holds for 96 of the 98 snarks of order $36$ with circular flow number $5$.  

It is important to stress that our method focuses on the presence of some structures which force the circular flow number of a snark to be large. Indeed, we describe a way to obtain graphs, not necessarily cubic, which are cyclically $4$-edge-connected graphs with circular flow number 5. Each such graph can then be transformed into a snark by a suitable expansion of some of its vertices (see Section~\ref{sec:knownmethods} and the Appendix for a precise description). 
The construction of such graphs is the main purpose of our method: all snarks obtained starting from them will have circular flow number 5, since the expansion of a vertex does not decrease the circular flow number (cf.\ Proposition~\ref{pro:expansion}). 
A paradigmatic case is that all 25 snarks of order 34 can be viewed as suitable expansions of the same unique graph. For the reader's convenience, an example taken from~\cite{EMT16} with a complete description of such a procedure is reported in the Appendix. 
To keep things concise, we will not specify every time how we can obtain a snark starting from a given graph. Moreover, along the entire paper, we only prove that our method produces graphs with circular flow number {\em at least} 5 as this is implicitly sufficient to prove that its circular flow number is {\em exactly} 5 if Tutte's Conjecture is indeed true. 

This paper is organised as follows.
In Section~$\ref{sect:generalised_edges}$, we present the terminology and notations introduced in~\cite{EMT16} that we are going to use in the rest of the paper. Moreover, we give a complete answer to Problem~7.3 from~\cite{EMT16}.
In Section~$\ref{sect:constructions}$, we summarise known constructions of graphs with circular flow number 5 and propose new ones. Afterwards, we prove that all such methods are nothing but particular instances of a more general construction that we introduce  into detail here.
Section~\ref{sect:main_results} and Section~\ref{sect:even_subg} are devoted to a complete analysis of many possible instances of the introduced method: a summary of all results obtained in these two sections is Theorem~\ref{theo:summary}, which is the main result of this paper.
Finally, in Section~\ref{sect:final_remarks}, we show the results of our computations. On one side, they confirm that our method is a good tool to produce several examples of snarks with circular flow number $5$, but, on the other hand, we find two snarks of order 36 which seem to not fit our description. This suggests that the variety of snarks with circular flow number 5 could be very large.

\section{Generalised edges and open $5$-capacity}\label{sect:generalised_edges}
	
This section is mainly devoted to a review of the main results and definition presented in~\cite{EMT16}, that will be needed later on. We refer to the same work by Esperet, Tarsi and the third author for a complete proof of the results of this section. Furthermore, we give a complete answer to Problem~7.3 in \cite{EMT16} which was left as an interesting open problem.

First of all we introduce the definition of modular flow and we recall that the existence of such a kind of flow is equivalent to the existence of a $r$-CNZF.

\begin{definition}
A \textbf{circular nowhere-zero modular r-flow}, or \textbf{r-MCNZF}, in a graph $G$, is an assignment $\phi \colon E \to [1,r-1] \subseteq \R/ r \Z$ together with an orientation of $G$, such that, for every $v\in V$, \[ \sum_{e\in E^{+}(v)} \phi(e) = \sum_{e\in E^{-}(v) }\phi(e) \mod{r}, \] where $E^{+}(v)$ and $E^{-}(v)$ denote the sets of ingoing and, respectively, outgoing arcs from a vertex $v$ in the selected orientation of $G$.
\end{definition}
		
\begin{proposition}\label{prop_flusso-intero-circolare}
An $r$-CNZF in a graph $G$ exists if and only if there exists an $r$-MCNZF. 
\end{proposition}
		
The following proposition gives an important tool that will be central in several proofs of the present paper.
		
\begin{proposition}\label{prop:circ_F_N_sub_Flow}
For a graph $G$, $\Phi_{c}(G)<r$ if and only if there exists an $r$-MCNZF $\phi$ in $G$ such that $\phi\colon E \to (1,r-1)$.
\end{proposition}
		
Following the notation used in~\cite{EMT16},  a flow which satisfies the condition in Proposition~\ref{prop:circ_F_N_sub_Flow} will be called a \emph{sub-$r$-MCNZF}. 

Let $r\in \mathbb{R}$ and consider $\mathbb{R}/r\mathbb{Z}$, the group of real numbers modulo $r$. This is commonly represented by a circle of length $r$, where $r$ coincides with $0$, and an open interval $(a,b)\subseteq \R / r \Z$ denotes the set of numbers covered when traversing clockwise from $a$ to $b$, with $a,b$ not included; closed intervals are denoted in a similar way. In particular $(x,x)$ is defined to be $\R/r \Z - \{x\}.$
		
We will focus on the case $\R / 5 \Z$.
The set of all \emph{integer open intervals} of $\R/ 5 \Z$, i.e.\ all intervals $(a,b)$ where $a,b\in \Z$, is denoted by $I_{5} := \{ (a,b)\subseteq \R/5 \Z\colon a,b\in \Z \}$.
A subset $X\subseteq \R / 5 \Z$ is called \emph{symmetric} if and only if $x\in X \iff -x \in X. $
We denote by $SI_{5}$ the family of all symmetric subsets of $\R/ 5 \Z$ which can be obtained as union of elements of $I_{5}$, that is:
$$SI_{5}:= \{ I \subseteq \R/ 5 \Z\colon I \text{ is symmetric and } I = \cup (a,b), (a,b) \in I_{5} \}.$$
		
\begin{definition}
The {\bf measure} of $A\in SI_{5}$, denoted by $\Me(A)$, is the number of unit intervals contained in $A$.
\end{definition}
		
The description of all constructions in the next sections of the present paper will be given making use of the definition of \textit{generalised edge}.
		
		\begin{definition}
			A {\bf generalised edge}, $G_{xy}$ is a graph $G$ with at least two vertices together with a pair of distinct vertices $x,y\in V(G)$ called its \bf{terminals}.
		\end{definition}
			
		Consider the generalised edge $G_{xy}$ and define a new graph (in general it can have multiple edges) $G_{xy}^{+}$ by adding a new edge $e^{+}=xy$ to $G_{xy}$.
		
		\begin{definition}
			The {\bf open $5$-capacity} of $G_{xy}$ is \[ CP_5 (G_{xy}) := \{ f(e^{+})\colon f \text{ is a modulo } 5 \text{-flow in } G_{xy}^{+} \text{ and } f|_E\subseteq (1,4) \}.  \]
		\end{definition}
				
		$CP_5 (G_{xy})$ is actually the set of values in $\R/5 \Z$, that can ``pass through'' $G$ from the source terminal $x$ to the sink terminal $y$, under all possible orientations of $G$, requiring that the flow capacity of every edge is restricted to the open interval $(1,4)$.
		
		From now on, since we are going to deal only with the case of modular $5$-flows, we will for simplicity refer to the open $5$-capacity of a generalised edge just as the capacity of that generalised edge.
		
		A strong relation between the concept of capacity of a generalised edge and the set $SI_5$ is given in the following lemma (see~\cite{EMT16}):
				
		\begin{lemma}\label{lem:subsetSI5}
			If $G_{xy}$ is a generalised edge, then $CP_{5}(G_{xy})\in SI_{5}$.
		\end{lemma}
		
		In view of the previous lemma, a generalised edge $G_{xy}$ having $5$-capacity $A \in SI_5$ is said to be an $A$-edge.  
		A standard edge (the graph with two vertices and one edge) is then a $(1,4)$-edge, but there exist infinitely many $(1,4)$-edges which are not isomorphic to it.  
		
		It is clear that any graph $G$ can be viewed as a union of generalised edges having disjoint vertex-sets except, possibly, for their terminals. Also note that the same graph could admit several different representations with different sets of generalised edges: a trivial representation is obtained by considering every edge of $G$ as a $(1,4)$-edge; on the opposite side, we can consider the entire graph $G$ and any two of its vertices as a generalised edge itself.   
		
		\begin{definition}
			Consider a pair $(H, \sigma)$, where $H=(V(H),E(H))$ is a graph and $\sigma: E(H)\to SI_5$ is a map that associates to each edge $uv\in E(H)$ a subset $\sigma(uv)\in SI_5$. 
			We denote by $H^\sigma$ the family of all possible graphs which can be obtained by replacing every edge $uv$ of $H$ with a $\sigma(uv)$-edge with terminals $u$ and $v$.
			We will refer to such a $\sigma$ as the { \bf capacity function} defined on $H$. 
		\end{definition}
		
		\begin{remark}
		Every graph $G$ belongs to the family $G^\sigma$ where $\sigma$ is the constant capacity function: $\sigma(e)=(1,4)$ for each $e \in E(G)$.
		\end{remark}
		
		The following proposition will play a crucial role in what follows. 
		
		\begin{proposition}\label{pro:Hsigma}
		A graph $G \in H^\sigma$ admits a sub-$5$-MCNZF if and only if $H$ admits a flow $f$ such that $f(e) \in \sigma(e)$, for all $e \in E(H)$.
		\end{proposition}
		
		The previous proposition also says that if a graph $G$ in $H^\sigma$ has circular flow number $5$, then all graphs in $H^\sigma$ have the same property. Hence, in order to find graphs with circular flow number 5, we can work on the pairs $(H,\sigma)$ instead of working on each specific graph $G$ of the family.   
		Mainly for this reason, we will make use of the following definition in the rest of the paper.
		
		\begin{definition}
		Given the family $H^\sigma$, we call a {\bf $\sigma$-faithful flow} in $H$ any flow in $H$ which satisfies the condition of Proposition~\ref{pro:Hsigma}.   
		\end{definition}

		\subsection{Every element of $SI_5$ is graphic}
		
		Lemma~\ref{lem:subsetSI5} shows that the open capacity of a generalised edge is an element of $SI_5$. 
		One of the open problems proposed in~\cite{EMT16} (i.e.\ Problem~7.3) is the determination of all elements of $SI_5$ which are the open $5$-capacity of a generalised edge. In order to study such a problem the following definition naturally arises:
	
		\begin{definition}
			Let $A\in SI_{5}$. If there exists an $A$-edge, then $A$ is called {\bf graphic}. We denote by $GI_{5}\subseteq SI_{5}$ the set of all graphic elements of $\R/5 \Z.$
		\end{definition}
		
		Two operations to produce new elements of $GI_5$ starting from the known ones are presented in~\cite{EMT16} which are used to prove the following proposition.
		
		\begin{proposition}\label{GI_k_algebra}
		$GI_{5}$ is a closed subfamily of $SI_{5}$ with respect to sum and intersection.
		\end{proposition}
		
		Only $5$ sets in $SI_{5}$ were not proved to be graphic in~\cite{EMT16}, more specifically those obtained by removing the two elements $\{2,3\}$ from the sets of $SI_{5}$ containing them.

		Now, we completely answer the question posed in~\cite{EMT16} by showing that also the remaining five sets of $SI_5$ are graphic, that is $GI_5=SI_5$. 
		We will make use of the following remark:
		\begin{remark}\label{rem:<5}
			Let $A\in GI_k$ and $H$ be an $A$-edge. Then $0\in CP_5(H)$ if and only if $\Phi_c(H)<5$.
		\end{remark}
		
		\begin{theorem}
			$GI_5 = SI_5 $.
		\end{theorem}
		\begin{proof}
		Consider the generalised edge $G_{uv}$ such that $G^+_{uv}$ is the complete graph with four vertices, and denote by $s$ and $t$ the other two vertices of $G$. 
		It is sufficient to prove that $CP_5(G_{uv})=\R/5\Z-\{ 2,3 \}$: indeed this would mean that $\R/5\Z-\{ 2,3 \}\in GI_5$ and, since $GI_5$ is closed under intersection, all remaining intervals could be conveniently generated. 
		First of all, let us show that $0$, $1$ and $2.5$ are elements of $CP_5(G_{uv})$.
		Since $\Phi_c(G)<5$, $0\in CP_5(G_{uv})$ follows from Remark~\ref{rem:<5}. For a sufficiently small $\epsilon>0$, we explicitly construct two flows in Figure~\ref{GI5} such that the flow value of $uv$ is $1$ (on the left) or $2.5$ (on the right).
				
		\begin{figure}[htb!]
			 \centering
			 \includegraphics[width=10cm]{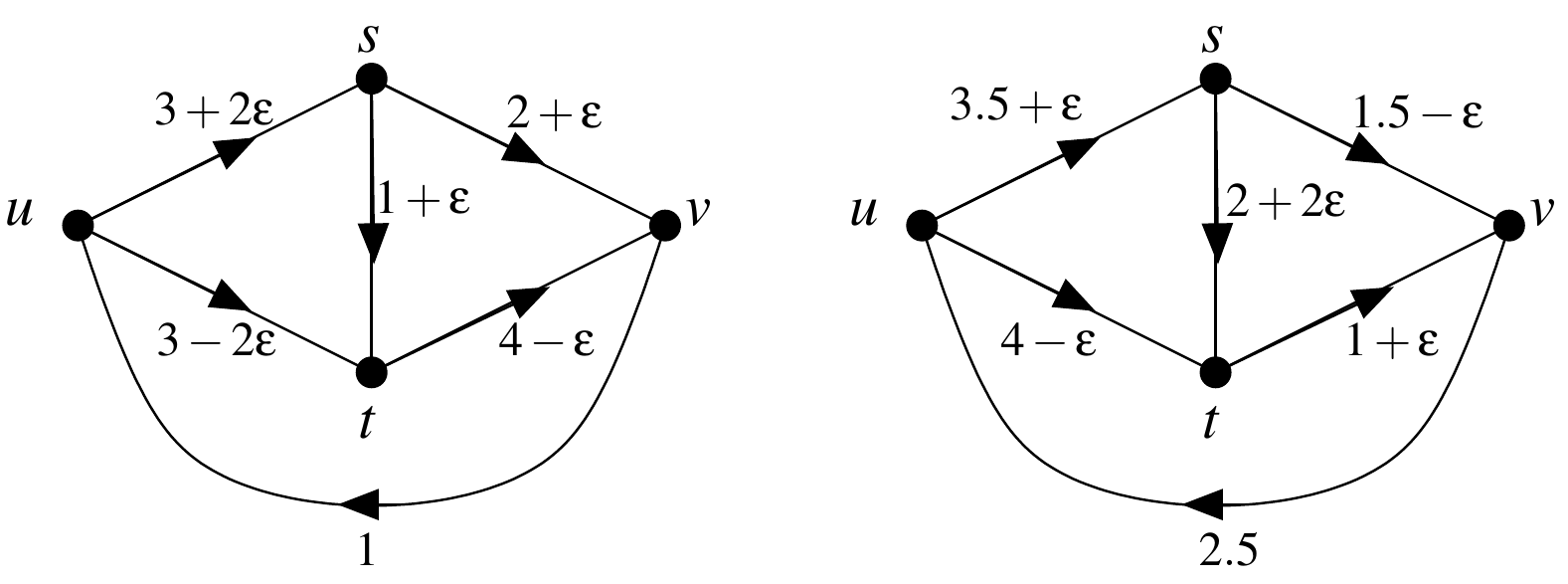}
			 \caption{Two flows of $K_4$ with prescribed flow  values on the edge $uv$.}\label{GI5}
		\end{figure}
			
		Hence, thanks to the openness and symmetry of the open capacity, we have proved that $\R/5\Z-\{ 2,3 \} \subseteq CP_5(G_{uv})$. In order to prove our assertion we need to show that $2\notin CP_5(G_{uv})$ (and by symmetry we also obtain $3\notin CP_5(G_{uv})$). 
		Take the same orientation of the edges of $G^+_{uv}$ shown in Figure~\ref{GI5} and suppose, by contradiction, that there is a flow $\phi$ in $G^+_{uv}$ such that $\phi|_E \subseteq (1,4)$ and $\phi(uv)=2$.
			
		From the relation $2= \phi(uv)= \phi(us)+\phi(ut)$, it follows that both $\phi(us),\phi(ut)\in (3,4)$. Similarly from $2= \phi(uv)= \phi(vs)+\phi(vt)$, we deduce that both $\phi(vs),\phi(vt)\in (3,4)$. Then, since $\phi(us)=\phi(st)+\phi(sv)$, it follows that $\phi(st)\in(4,1)$ a contradiction.
		\end{proof} 		

		The main result of this section says that every element of $SI_5$ is the $5$-capacity of a suitable generalised edge. Hence, from now on, every time we will consider a pair $(H,\sigma)$, we have no general restriction on the values assumed by $\sigma$.

		\section{Construction of graphs with circular flow number~5}\label{sect:constructions}
		
		We denote by $F_{\ge5}$ the family of graphs with circular flow number greater than or equal to $5$, and by $S_{\ge5}$ the subfamily of $F_{\ge5}$ consisting of snarks.
		
		\subsection{Known methods}\label{sec:knownmethods}
		
		Some of the constructions of snarks in $S_{\ge5}$ presented in~\cite{EMT16} make use of the following lemma. We report it into detail since it will be used  in the next section to describe some new methods.
		
		\begin{lemma}\label{Lemma_Cammino_Alternato}
			Consider a pair $(H,\sigma)$, where $H$ is a graph and $\sigma$ a capacity function defined on $H$. Suppose that $P$ is a path in $H$ with $\sigma(e)=A$ for every $e \in E(P)$ and $\Me(A)=2$. Assume also that each internal vertex $v_i$ of $P$ is adjacent to exactly one vertex $v'_i$ of $H$ not in $P$. Finally, assume $\sigma(v_iv'_i)\subseteq (1,4)$ for every vertex $v_i$. If $P$ is a directed path in a suitable orientation of $H$ and $f$ is a flow in $H$ such that $f(e)\in \sigma(e)$ for every $e\in E(H)$, then $f$ assigns to adjacent edges of $P$ two values which lie in the two different unit intervals of $A$.
		\end{lemma}
		
		The main method presented in~\cite{EMT16} to produce graphs in $F_{\ge5}$ is a direct application of the following corollary. Also note that the method previously presented in~\cite{MaRa06} is a particular case of the same corollary.
		
		\begin{corollary}\label{Cor:P_odd_cycle_1}
			Consider a pair $(H,\sigma)$, where $H$ is a graph and $\sigma$ a capacity function defined on $H$. Suppose that $C$ is an odd cycle in $H$ with $\sigma(e)=A$ for every $e \in E(C)$ and $\Me(A)=2$. Assume also that each vertex $v_i$ of $C$ is adjacent to exactly one vertex $v'_i$ of $H$ not in $C$. Finally,  assume $\sigma(v_iv'_i)\subseteq (1,4)$ for every vertex $v_i$. Then, $\Phi_{c}(G)\ge 5$ for all $G \in H^\sigma$.
		\end{corollary}
		In particular, if $H$ has sufficiently large girth and connectivity, we can construct an element of $S_{\ge5}$ starting from a suitable graph $G \in H^\sigma$. As already remarked, the standard trick to obtain a snark starting from an element $G$ of $H^\sigma$ is by applying an expansion operation. 
		
		\begin{definition}
			An {\bf expansion} of a vertex $x \in G$ into a graph $K$ is obtained by the replacement of $x$ by $K$ and by adding as many edges  with one end in $V(G-x)$ and the other end in $V(K)$ as the degree of $x$ in $V(G)$.
		\end{definition}
		
		Many different expansions can be performed on the same graph $G$, but it is well-known that these expansions do not decrease the circular flow number.
		
		\begin{proposition}\label{pro:expansion}
			Let $G'$ be a graph obtained with an expansion of a vertex of $G$. Then, $\Phi_{c}(G')\ge \Phi_{c}(G).$
		\end{proposition}
		
		An operation that trivially preserves the circular flow number is defined as follows:
		
		\begin{definition}
			If $x\in G$ is a vertex of degree $2$, let $N_G(x)=\{y,z\}$. Then {\bf smoothing} $x$ means removing the vertex $x$ and adding to $G$ the new edge $yz$.
		\end{definition}
		
		In the Appendix we give an example which shows how previous operations permit to obtain a snark with circular flow number 5.
		
\subsection{New methods}
	
		In previous section we have briefly described a method from~\cite{EMT16} to generate new graphs in $F_{\ge5}$. That construction gives ways to generate new members of $F_{\ge5}$ starting from graphs having particular subgraphs and capacity functions.
	
		Our next goal is to present some new constructions and, after that, to suggest a unified description of all methods described in this section and in the previous one. The final aim is a significant reduction of redundancy and a much better control on the graphs that can be generated. Such a unified description will be studied into detail in the next section and it is the main goal of this paper.
		
		Let us begin with new corollaries from Lemma~\ref{Lemma_Cammino_Alternato}. Each of them produces a new method to generate elements of $F_{\ge5}$.
		
		\begin{corollary}\label{cor_costruction_2}
			Consider a pair $(H,\sigma)$, where $H$ is a graph and $\sigma$ a capacity function defined on $H$. Suppose that $P$ is a path in $H$ with $\sigma(e)=(4,1)$ for every $e \in E(P)$. Also assume that each internal vertex $v_i$ of $P$ has degree $3$ in $H$ and that $\sigma(e)\subseteq (1,4)$ for every edge $e$ not in $E(P)$ and incident to an internal vertex of $P$. 
			If two internal vertices of $P$ at even distance on $P$ are adjacent, i.e.\ there exists an edge $e$ not in $P$ connecting two internal vertices of $P$ and forming  an odd cycle with the edges of $P$,  then $\Phi_{c}(G)\ge 5$ for all $G \in H^\sigma$.
		\end{corollary}

		\begin{proof}
			Give an orientation to $P$ in such a way that it becomes a directed path in $H$ and let $f$ be a flow in $H$ such that $f(e)\in \sigma(e)$ for each $e\in E(H)$. Then, by Lemma~\ref{Lemma_Cammino_Alternato}, the edges of $P$ take values alternately from the two unit intervals $(4,0)$ and $(0,1)$: a contradiction arises from the fact that $f(e)$ must stay at the same time in $(4,0)-(0,1)=(3,0)$ and in $(0,1)-(4,0)=(0,2)$. Hence such a flow $f$ cannot exist and so, by Proposition~\ref{pro:Hsigma}, every $G\in H^{\sigma}$ cannot have a sub-$5$-MCNZF. 
		\end{proof}

		\begin{corollary}\label{cor_costruction_3}
		Consider a pair $(H,\sigma)$, where $H$ is a graph and $\sigma$ a capacity function defined on $H$. Suppose that $P_1$ and $P_2$ are distinct paths in $H$ with $\sigma(e)=(4,1)$ for every $e \in E(P_1)\cup E(P_2)$. Also assume that each internal vertex $v_i$ of these paths has degree $3$ in $H$ and $\sigma(e)\subseteq (1,4)$ for every edge $e$ not in $E(P_i)$ and incident to an internal vertex of $P_i$. 
		
		If two internal vertices of $P_1$ at even distance (on $P_1$) are adjacent, respectively, to two internal vertices of $P_2$ at odd distance (on $P_2$), i.e.\ there exist two edges not in $P_1 \cup P_2$ connecting two internal vertices of $P_1$ to two internal vertices of $P_2$ and forming  an odd cycle with some edges of $P_1$ and $P_2$,  then $\Phi_{c}(G)\ge 5$ for all $G \in H^\sigma$.
		\end{corollary}
		
		\begin{proof}
			Give a suitable orientation to $H$ that makes both $P_1$ and $P_2$ directed paths and suppose that there exists a flow $f$ in $H$ such that $f(e)\in \sigma(e)$ for each $e\in E(H)$. We can assume without loss of generality that $P_1= x_0\dots x_s$, $P_2 = y_0\dots y_{t}$, with $t>s$, and that the two edges in the hypothesis are $x_1y_1$ and $x_{s-1}y_{t-1}$. By Lemma~\ref{Lemma_Cammino_Alternato}, the edges of each $P_j$ take values alternately from $(4,0)$ and $(0,1)$, but the presence of the edge $x_1y_1$ between them obliges the paths to start with different intervals, i.e.\ $x_0x_1\in (0,1)$ (resp. $(4,0)$) if and only if $y_0y_1\in (4,0)$ (resp. (0,1)). Since $s-1$ and $t-1$ have different parity, the values $x_{s-2}x_{s-1}-x_{s-1}x_s$ and $y_{t-2}y_{t-1} - y_{t-1}y_{t}$ belong to the same unit interval $(0,1)$ or $(4,0)$, whence there is no orientation of $x_{s-1}y_{t-1}$ such that the flow $f$ does exist. Therefore, by Proposition~\ref{pro:Hsigma}, every $G\in H^{\sigma}$ cannot have a sub-$5$-MCNZF.
		\end{proof}
		
		We conclude this section by proving that all previous results can be slightly generalised when we replace $(4,1)$ with any of its subsets. This fact is an obvious consequence of the following more general proposition.
		
		\begin{proposition}\label{pro:A'_subset_A}
			Let $H$ be a graph and let $\sigma_1$ and $\sigma_2$ be two capacity functions defined on $H$. Assume that $\sigma_2(e)\subseteq \sigma_1(e)$ for all $e \in E(H)$. If $\Phi_c(G)\ge 5$ for $G \in H^{\sigma_1}$, then $\Phi_c(G')\ge 5$ for $G' \in H^{\sigma_2}$.
		\end{proposition}
		\begin{proof}
		It is sufficient to notice that $\sigma_2$ is a more restrictive capacity function. Then, starting from a sub-5-MCNZF of $G'$, we can reconstruct a sub-5-MCNZF of $G$, a contradiction. 
		\end{proof}

		In particular, thanks to Proposition~\ref{pro:A'_subset_A}, both corollaries presented in this section hold as well if some $(4,1)$-edges are replaced with $(4,0)\cup (0,1)$-edges.
		
		\subsection{A unified description}
		
		Now we suggest a possible unified description of the three methods arising from Corollary~\ref{Cor:P_odd_cycle_1}, Corollary~\ref{cor_costruction_2} and Corollary~\ref{cor_costruction_3}.
		
		Our approach is the following: we consider the subgraph of $(H,\sigma)$ described in the corresponding corollary which forces all graphs in $H^\sigma$ to have circular flow number at least 5, and we contract the remaining part of $H$ into a unique vertex.

		\begin{figure}[htb!]
			\centering
			\includegraphics[scale=0.5]{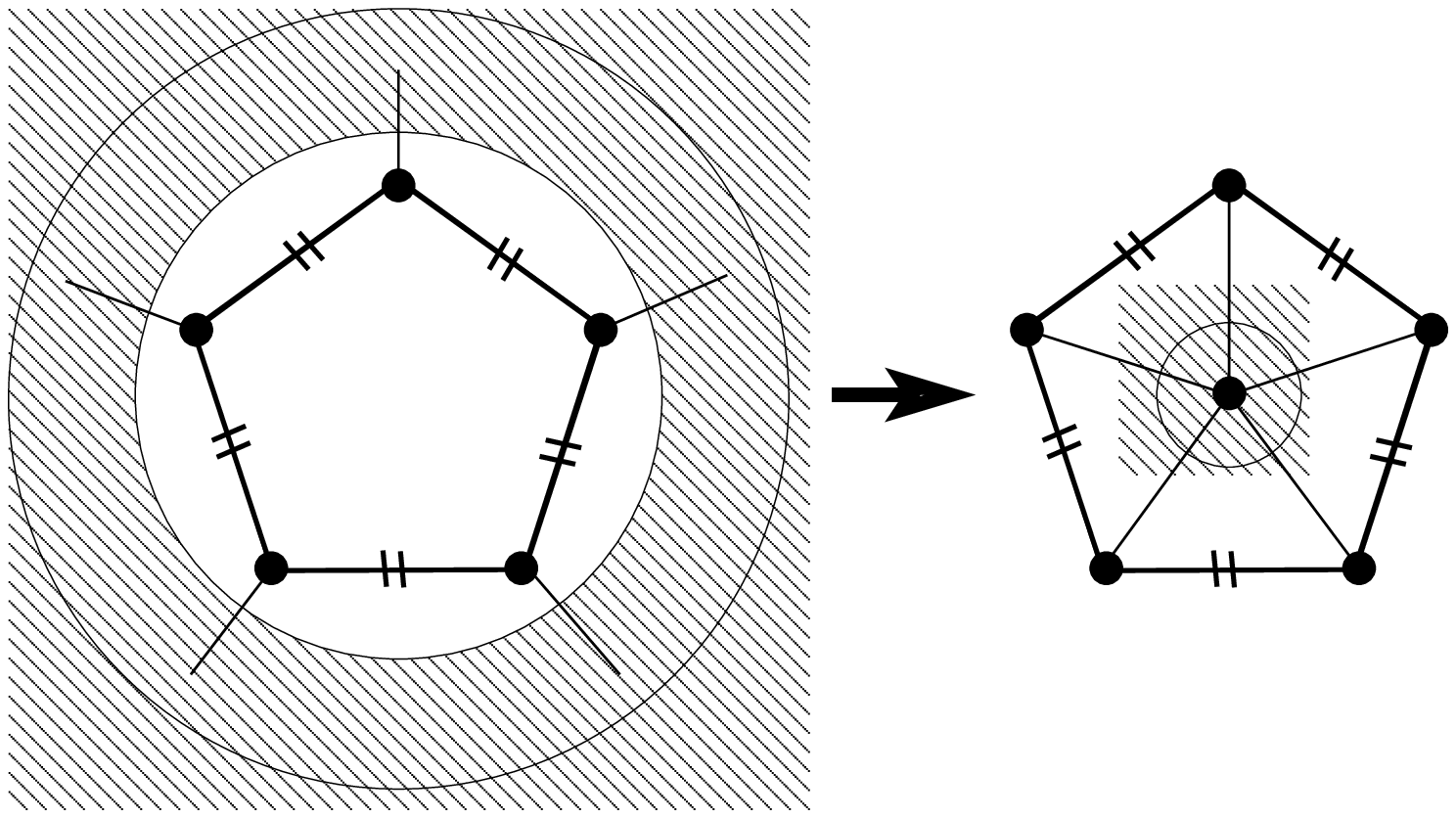}
			\caption{Example of a reduction based on Corollary~\ref{Cor:P_odd_cycle_1}.}\label{red_1}
			\includegraphics[scale=0.5]{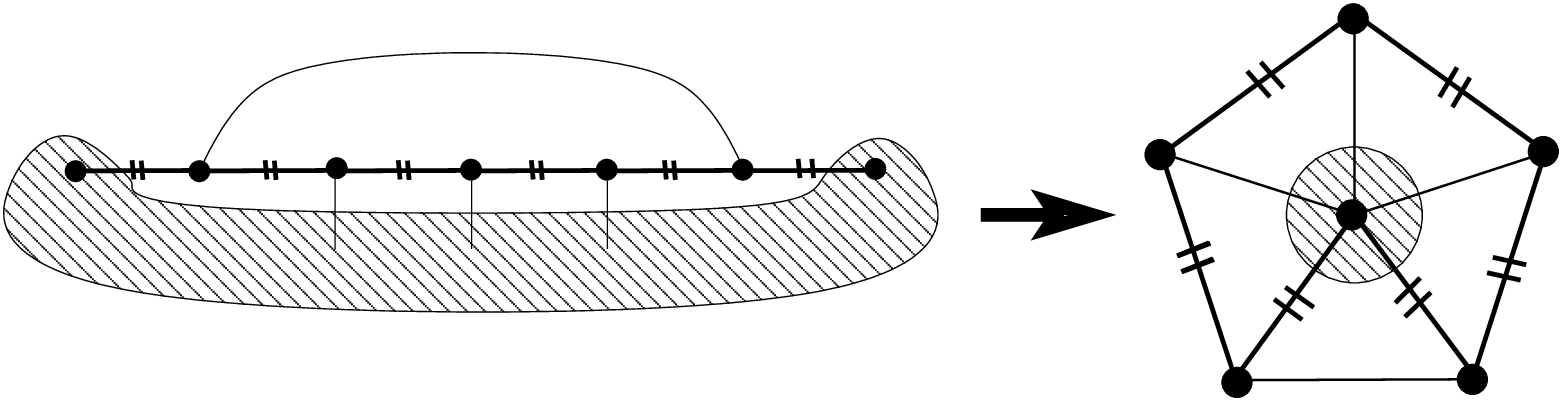}
			\caption{Example of a reduction based on Corollary~\ref{cor_costruction_2}.}\label{red_2}
			\includegraphics[scale=0.5]{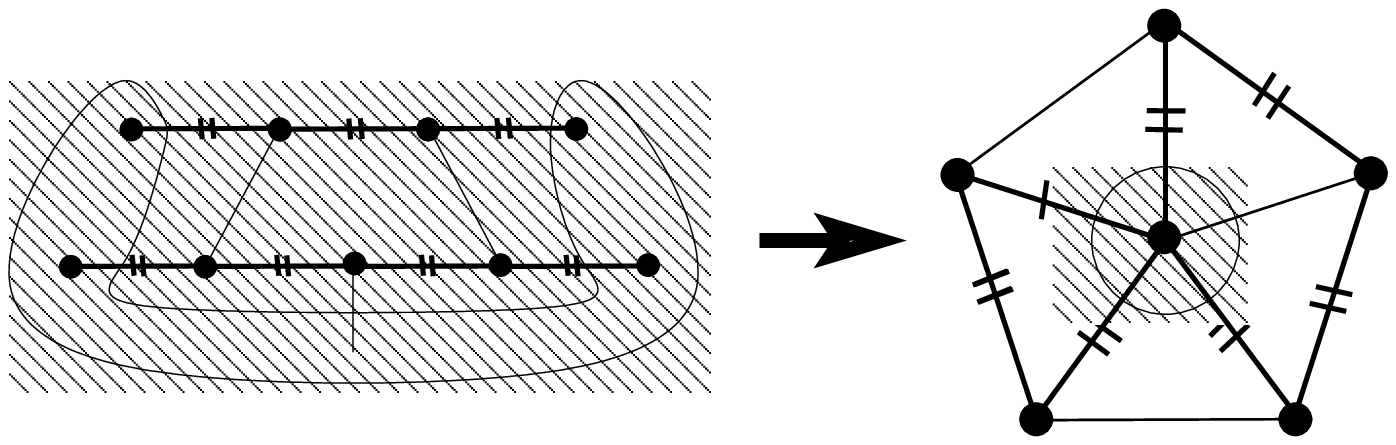}
			\caption{Example of a reduction based on Corollary~\ref{cor_costruction_3}.}\label{red_3}
		\end{figure}
			
		The key point is that the resulting graph is a wheel in all previous reductions (see Figures~\ref{red_1}-\ref{red_3}). 
		
		Moreover, the subgraph induced by the edges with capacity contained in $(4,1)$ (denoted by a double signed edge in the figures) is an even subgraph of the wheel, where an \emph{even subgraph} of a graph $H$ is a subgraph where all vertices have even degree.

		More precisely, we can summarise all previous reductions in the following list:

		\begin{itemize}
			\item Corollary~\ref{Cor:P_odd_cycle_1}: in case $C$ is a $(2n+1)$-cycle, this generates graphs which belong to $W^\sigma_{2n+1}$, where $\sigma(e) \subseteq (4,1)$ for all $e$ of the external cycle of $W_{2n+1}$ and $\sigma(e) \subseteq (1,4)$ for all other edges of the wheel (see Figure~\ref{red_1});
			\item Corollary~\ref{cor_costruction_2}: in case $P$ is a path with $2n+3$ vertices, this generates graphs which belong to $W^\sigma_{2n+1}$, where $\sigma(e) \subseteq (4,1)$ for all $e$ of a Hamiltonian cycle of $W_{2n+1}$ and $\sigma(e) \subseteq (1,4)$ for all other edges of the wheel (see Figure~\ref{red_2});
			\item Corollary~\ref{cor_costruction_3}: in case $P_1$ and $P_2$ together have $(2n+5)$ vertices, this generates graphs which belong to $W^\sigma_{2n+1}$, where $\sigma(e) \subseteq (4,1)$ for all $e$ of a suitable even subgraph of $W_{2n+1}$ and $\sigma(e) \subseteq (1,4)$ for all other edges of the wheel (see Figure~\ref{red_3}).
		\end{itemize}
			
		The next section is devoted to an exhaustive analysis of all instances arising from the new approach we have introduced in this section.
		
		\section{Main results}\label{sect:main_results}
		
		In the previous section, we remark that all construction methods of cubic graphs having circular flow number at least 5 always produce a graph which is a suitable expansion of a graph in $W_{n}^\sigma$ where $\sigma$ is a capacity function which assigns a subset of $(4,1)$ to all edges of a given even subgraph of the wheel and $(1,4)$ to all other edges of $W_n$. 
		
		So it is very natural to ask in general whether, given a wheel $W_n$ and a capacity function $\sigma$ such that $\sigma(e)\subseteq (4,1)$ for all edges of an even subgraph $J$ of $W_n$ and $\sigma(e)=(1,4)$ otherwise, we obtain that a graph $G$ which belongs to $W_n^\sigma$ has circular flow number at least 5.
		
                \begin{remark}
                In this section, we consider only the case in which $\sigma(e)=(4,1)$ for all edges of the even subgraph $J$ and we use the notation $J_{(4,1)}$ to stress the fact that all edges of $J$ have capacity $(4,1)$. Anyway, it follows by Proposition~\ref{pro:A'_subset_A} that all results we are going to present also hold in the more general case $\sigma(e) \subseteq (4,1)$ for some edges of $J$.
		\end{remark}
		
		More precisely we can consider the following problem:
		
		\begin{problem}
		Given a wheel $W_n$ with $n+1$ vertices and $J$ a non-empty even subgraph of $W_n$, establish for each integer $n$ and each possible even subgraph $J$, if a graph $G \in W_n^\sigma$ has circular flow number at least 5, where $\sigma(e)=(4,1)$ if $e \in E(J)$ and $\sigma(e)=(1,4)$ otherwise. We will denote such a family of graphs by $(W_n,J_{(4,1)})$.  
		\end{problem}
		
		In what follows, we give a complete solution for all possible instances of this problem.
		
		We call a $(4,1)$-edge an edge of $W_n$ which belongs to $J$ and a $(1,4)$-edge an edge of $W_n$ which does not belong to $J$. 
		
\begin{figure}
			\centering
			\includegraphics[scale=0.6]{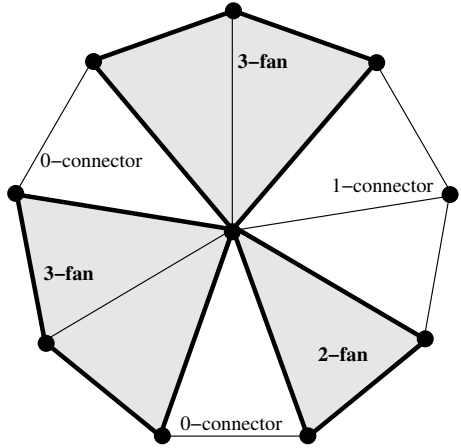}
			\caption{An example of fans and connectors in $W_9$, where the bold edges represent edges of $J$ and all others are $(1,4)$-edges.}\label{fanconn}
		\end{figure}

		In order to describe the structure of $J$, we introduce the following two useful definitions (see also Figure~\ref{fanconn}).
		
		\begin{definition}
			Given the family $(W_n,J_{(4,1)})$, for an integer $l\ge 2$, an {\bf $l$-fan} $F_l$ in $W_{n}$ is a subgraph induced by all vertices of an $(l+1)$-cycle consisting of edges of $J$ and passing through the central vertex $v_c$ of $W_n$.
			\end{definition}
		
		\begin{definition}
			Given the family $(W_n,J_{(4,1)})$, for an integer $m\ge0$, an {\bf $m$-connector} $C_m$ is a subgraph of $W_{n}$ induced by all the $m+1$ edges of a maximal path $P_{m+2}$ of $(1,4)$-edges of the external cycle of $W_n$ and all $(1,4)$-edges of type $uv_c$, where $u$ is a degree $2$ vertex of $P_{m+2}$ and $v_c$ is the center of $W_n$.
		\end{definition}
		
		It is clear that every even subgraph, except if $J$ is the empty graph or the external cycle of $W_n$, can be described as a sequence of fans and connectors in $W_n$. From now on, when we refer to the connector (fan) \emph{following} a fan (connector), we are implicitely considering the clockwise order on fans and connectors in  $W_n$.  
		
		In what follows, we also use the following terminology: 
		\begin{itemize}
		\item We refer to the longest cycle of an $m$-connector, for $m\geq2$, as the \emph{external cycle} and, accordingly, we call its edges \emph{external edges}. For all $m \geq 0$, we call \emph{internal edges} all edges of an $m$-connector incident to $v_c$ which are not external. Finally, we call \emph{lateral edges} of an $m$-connector the two edges (only one in the case of a $0$-connector) that are neither external nor internal.
		
		\item We refer to the cycle of edges of $J$ in an $l$-fan as the \emph{external cycle} of the fan and, accordingly, we call its edges \emph{external edges}. We will refer to the $(1,4)$-edges of an $l$-fan as its \emph{internal} edges. Finally, we will call the \emph{first edge of the fan} the unique external edge of $J$ which is incident to $v_c$ and to a lateral edge of the connector preceding the fan. While we call the \emph{last edge of the fan} the unique external edge of $J$ which is incident to $v_c$ and to a lateral edge  of the connector following the fan.
		\end{itemize}
		
		\subsection{Flows in fans and connectors} 
		
		Now we furnish the description of some flows defined in $l$-fans and $m$-connectors that will be largely used in the following proofs. 
		
		In what follows, take three values $x,y,z \in \R/5 \Z$ such that $x\in(1,2)$, $y \in (1,2)$, $z\in (4,0)$, $x+z\in (0,1)$ and $x+2y \in (1,4)$: it is an easy check that such values do exist. Moreover, note also that $2y \in (1,4)$ and $x+y \in (1,4)$.
		
		\underline{Flow $f^+$ in an $l$-fan ($l$ even):} consider an $l$-fan $F_l$ with $l\geq2$ even, we assign a clockwise orientation to its external cycle and we define $f^+$ such that it assigns flow value $z$ and $z+x$ alternately to the edges of the external cycle starting from $v_c$ and following the orientation. Moreover, $f^+$ assigns flow value $x$ to all internal edges of $F_l$. Now, we consider the unique possible orientation of the internal edges such that $f^+$ is a zero-sum flow in each vertex of degree $3$ in $F_l$. Indeed, note that if $l$ is even then $f^+$ is also a zero-sum flow in $v_c$.
		
		\underline{Flow $\tilde{f}^+$ in an $l$-fan ($l$ odd):} consider an $l$-fan $F_l$ with $l\geq3$ odd, the flow $\tilde{f}^+$ assigns the orientation and the flow values exactly as $f^+$, but note that if $l$ is odd then  the difference between the inner flow and outer flow of $\tilde{f}^+$ in $v_c$ is exactly $2x$. 
		
\begin{figure}[htb!]
			\centering
			\includegraphics[scale=0.9]{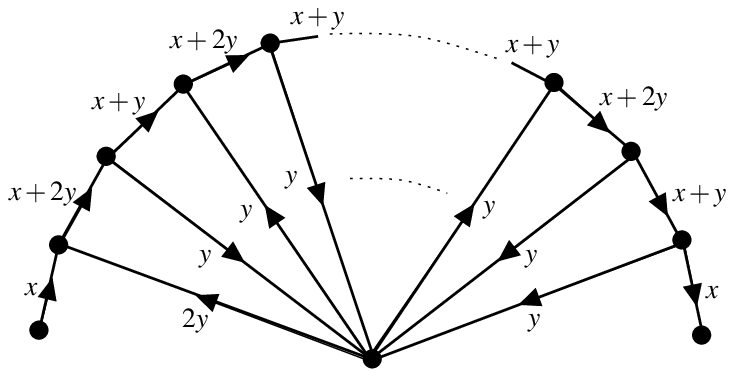} \hspace{1cm} \includegraphics[scale=0.9]{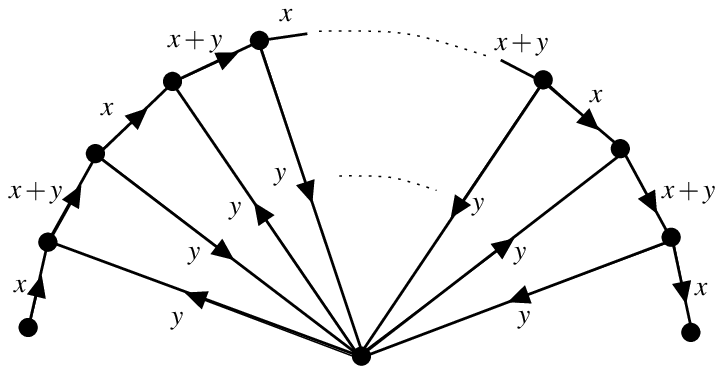}
			\caption{The flow $g^+$ in an odd connector (on the left) and in an even connector (on the right).}\label{Fig:g+}
		\end{figure}

		\underline{Flow $g^+$ in an $m$-connector (for $m \neq 1$):} consider an $m$-connector $C_m$. For $m=0$ we set the flow value of $g^+$ equal to $x$ on the unique edge of the connector, which is oriented clockwise in $W_n$. For $m>1$ we distinguish two cases according to the parity of $m$. 
		If $m>1$ odd, then we define $g^+$ as on the left of Figure~\ref{Fig:g+}, while if $m>0$ even then we define $g^+$ as on the right of Figure~\ref{Fig:g+}.

\begin{figure}[htb!]
			\centering
			\includegraphics[scale=0.9]{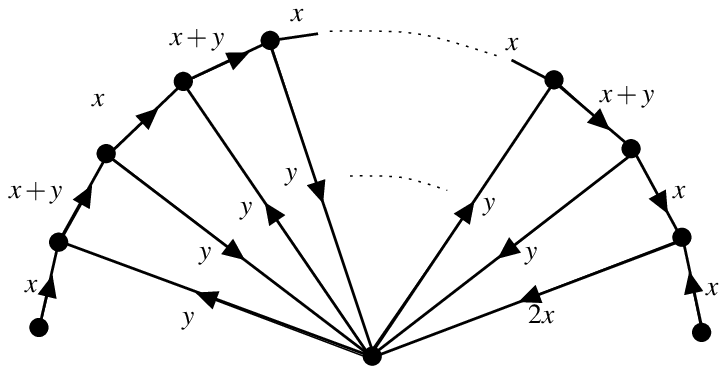} \hspace{1cm} \includegraphics[scale=0.9]{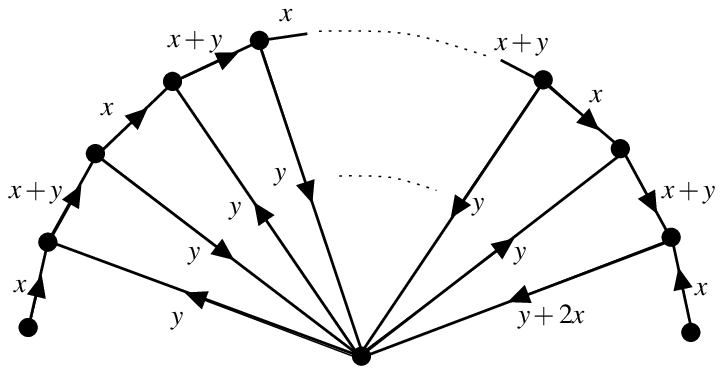}
			\caption{The flow $\tilde{g}^+$ in an odd connector (on the left) and in an even connector (on the right).}\label{Fig:gtilde+}
		\end{figure}

		\underline{Flow $\tilde{g}^+$ in an $m$-connector (for $m>0$):} consider an $m$-connector $C_m$.  
		For $m>1$ we distinguish two cases according to the parity of $m$. If $m>1$ odd, then we define $\tilde{g}^+$ as on the left of Figure~\ref{Fig:gtilde+}, while if $m>0$ even then we define $\tilde{g}^+$ as on the right of Figure~\ref{Fig:gtilde+}. In particular, note that for a $1$-connector, if we denote by $uv_c$ the unique edge incident to $v_c$, $\tilde{g}^+$ orients it towards $v_c$ and assigns it flow value $2x$. Moreover $\tilde{g}^+$ assigns to the other two edges flow value $x$ and orients them in such a way that $\tilde{g}^+$ is a zero-sum flow in $u$.

		\underline{Flows $f^-$, $g^-$ and $\tilde{g}^-$:} these flows are obtained from $f^+$, $g^+$ and $\tilde{g}^+$, respectively, by considering the same flow value on each edge of the corresponding flow and by reversing the orientation of each edge with respect to the orientation in the original one. 
		
		Consider the family $(W_n,J_{(4,1)})$ and let $f$ be a $\sigma$-faithful flow in $W_n$ which coincides with one of the flows $f^+$,$f^-$, $\tilde{f}^+$ or $\tilde{f}^+$ when we consider its restriction to an $l$-fan $F_l$ with $l<n$. 
		Consider the unique two lateral edges in $W_n$ which are incident to a vertex of the $l$-fan and assume that $f$ has flow value $x$ on both of them. Then, there is a unique way to orient these lateral edges, say $e_1$ and $e_2$, in such a way that $f$ is a zero-sum flow in all vertices distinct from $v_c$ of $F_l$. 
		If we have the flow $f^+$ ($f^-$) on the fan, both $e_1$ and $e_2$ have a clockwise (anticlockwise) orientation, otherwise in $\tilde{f}^+$ ($\tilde{f}^-$), they both point towards (away from) the fan.
		
		\begin{figure}[htb!]
			\centering
			\includegraphics[trim={0 0 3cm 0},clip,scale=0.6]{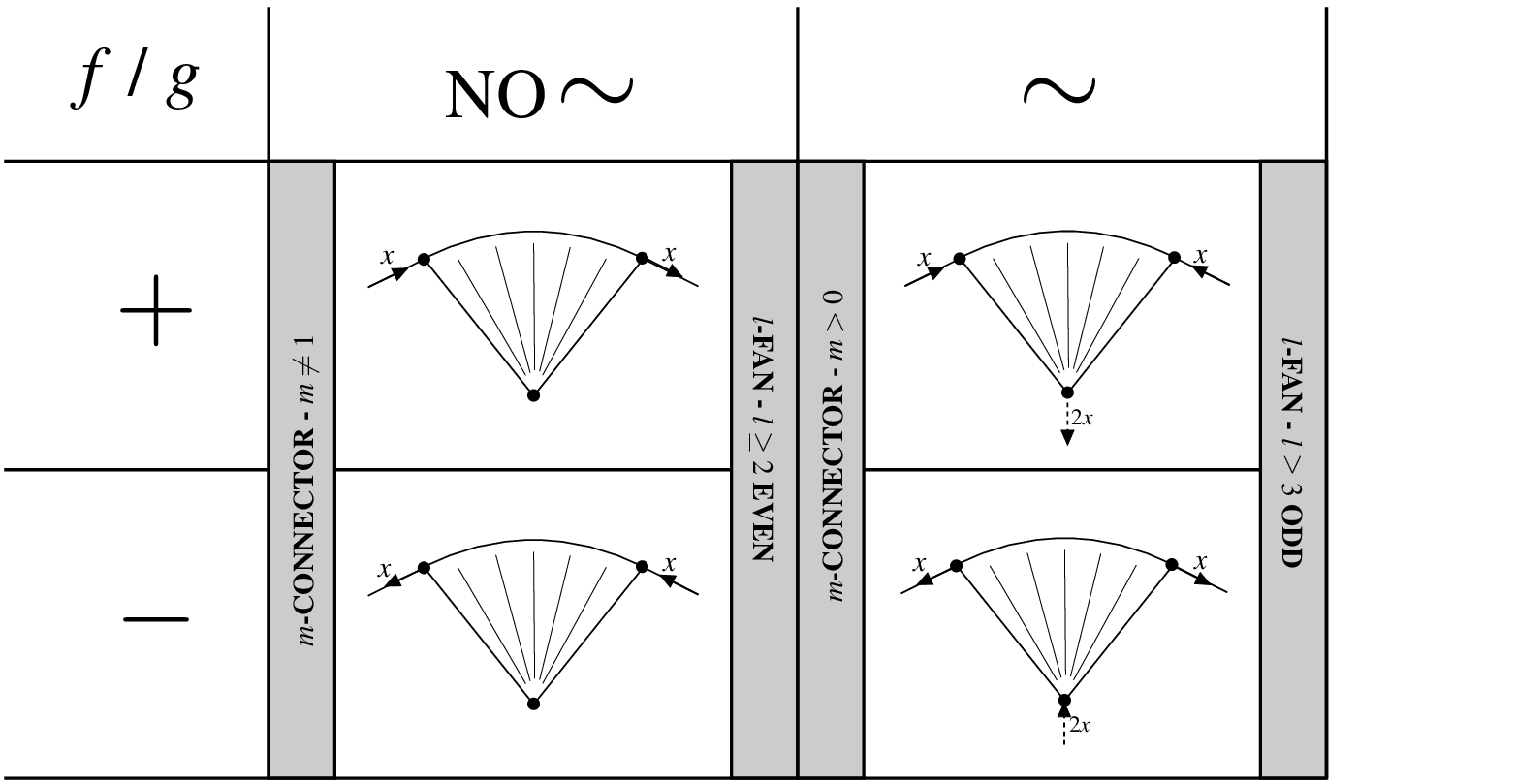}
			\caption{A summary of all flows defined in this section.}\label{Fig:flows}
		\end{figure}
		
		More in general, our notation for all previous flows is consistent with the following scheme from Figure~\ref{Fig:flows} which uses the following notation:
		\begin{itemize}
		\item The letters $f$ and $g$ denote flows on fans and connectors, respectively.
		\item The symbols $+$ and $-$  mean that the edges of the external cycles of $W_n$ are oriented in clockwise and anticlockwise direction in the corresponding flow, respectively.
		\item The symbol $\sim$ means that, in the vertex $v_c$, the absolute difference of the inner flow and outer flow in the corresponding subgraph is $2x$. More precisely, it is $2x$ in the case of all flows with symbol $+$, and $-2x$ in the case of all flows with symbol $-$. Whereas, the absence of $\sim$ means that the corresponding flow is a zero-sum flow in $v_c$.
		\end{itemize}
		
		In other words, in all flows with the symbol $\sim$ the two lateral edges of the connector (or the two lateral edges adjacent to the fan) have the same value $x$ and opposite orientation in the cycle of $W_n$, while in all flows without the symbol $\sim$ we have the same value $x$ on both edges and the same orientation in the cycle of $W_n$.
		
		We can briefly say that in the former case the flow \emph{reverses} the orientation, while in the latter case the flow \emph{preserves} the orientation.
				
		\subsection{Even subgraphs with edges of capacity $(4,1)$}
		
		In this section, we completely characterise the families $(W_n,J_{(4,1)})$ whose elements are graphs with circular flow number at least 5.
		When we speak about the sequence of fans and connectors given by the choice of $J$ in $W_n$ in the proofs, we will use the term \emph{component} to speak indifferently about either a fan or a connector of the decomposition. 
		
		\begin{proposition}\label{prop:subflow_even_wheel_41_edges}			
			If $G$ belongs to $(W_{n},J_{(4,1)})$ with $n$ even, then $\Phi_c(G)<5$.
		\end{proposition}
		
		\begin{proof} Set $n=2k$. If $J= W_{2k}-v_c$, we construct a $\sigma$-faithful flow $\phi$ in $W_{2k}$ and the assertion follows by Proposition~\ref{prop:circ_F_N_sub_Flow} and Proposition~\ref{pro:Hsigma}. Assign alternately the flow values $z$ and $x+z$ to all edges of the external cycle of $W_{2k}$ oriented clockwise and the flow value $x$ to all other edges with the unique orientation that makes $\phi$ a zero-sum flow at each vertex. Hence, we can assume without loss of generality that $v_c\in J.$					
		As already remarked, we can describe $J$ as a sequence of fans and connectors. Now, we recursively assign the flow on the components of such a decomposition. Firstly, select an $l$-fan and assign it the flow $f^+$ if $l$ is even and the flow $\tilde{f}^+$ if $l$ is odd. At each further step consider the following component and assign it a flow with these rules:
			\begin{itemize}
			\item if the component is a fan then assign a flow $f$, if it is a connector then assign a flow $g$; 
			\item if $l$ (or $m$) is odd, then assign a flow with $\sim$;  
			\item if the previous component has a flow with $\sim$, then assign a flow with opposite sign with respect to the flow in the previous component, otherwise a flow with the same sign. 
			\end{itemize}

			Since the wheel is even, the number of odd components is even. Hence, from our construction it follows that there is an even number of components with a flow that reverses the orientation. Hence the flows defined on each subgraph, altogether, induce a zero-sum flow on each vertex of the external cycle, and then also at the vertex $v_c$. Since we have defined a $\sigma$-faithful flow in $W_{2k}$, the assertion follows.
		\end{proof}
		
		Now we complete the characterization of even subgraphs $J$ such that a graph in $(W_n, J_{(4,1)})$, $n$ odd, has circular flow number at least $5$.
		
		\begin{theorem}\label{teo:odd_wheels_41_edges}			
			Let $G$ be a graph in $(W_{n},J_{(4,1)})$, with $n$ odd. Then \[ \Phi_c (G)\ge 5 \text{ if and only if } J \text{ has no } k \text{-connectors, for } k\ge2. \]
		\end{theorem}
		
		\begin{proof}			
			Assume that $J$ can be described by using only $0$-connectors and $1$-connectors.
			If we prove that a $\sigma$-faithful flow cannot be defined for a given orientation, then it cannot be defined for an arbitrary orientation. We orient all edges of the external cycle of $W_n$ and all edges of the external cycle of each fan in a clockwise way. Take an arbitrary orientation of all remaining edges and assume, by contradiction, the existence of a $\sigma$-faithful flow.
			First of all, recall that the flow value of the edges in the external cycle of a fan must be in $(4,1)$. Moreover, by Lemma~\ref{Lemma_Cammino_Alternato}, the flow value on two consecutive external edges of a fan must be in the two disjoint unit intervals $(4,0)$ and $(0,1)$ (except possibly for the two edges sharing the central vertex $v_c$). It follows that every lateral edge must have a flow  value either in the unit interval $(0,1)-(4,0)=(1,2)$ or in the unit interval $(4,0)-(0,1)=(3,4)$, because it is incident to two consecutive external edges of a fan.
			Note that, due to the chosen orientation,  if $l$ is even, lateral edges of an $l$-fan take values in the same unit interval, whereas, if $l$ is odd, they take values in each one of the two different unit intervals, respectively. 
			Analogously, consider the two lateral edges of a $1$-connector. It is easy to verify that these two edges must have flow values in $(1,2)$ and $(3,4)$, respectively. 
			Obviously, the unique lateral edge of a $0$-connector takes a flow value in one of those two unit intervals.
			Now consider an arbitrary lateral edge $e$, without loss of generality we can assume it has a flow value in $(1,2)$. By previous considerations, starting from $e$ we can establish in which interval between $(1,2)$ and $(3,4)$ each lateral edge lies: if the next component is either an odd fan or  a $1$-connector, then the flow value of the next lateral edge is in the other unit interval, otherwise it is in the same unit interval.
			Since in this case the number of odd components is odd, the flow value of $e$ should belong both to $(1,2)$ and $(3,4)$, a contradiction. 
			
			For the necessity, assume that there exists an $m$-connector $C_m$, with $m\geq2$. We show that there is a $\sigma$-faithful flow in $W_n$.
			We assign orientations and flow values to all fans and connectors, except for $C_m$, following exactly the same rules described in Proposition~\ref{prop:subflow_even_wheel_41_edges}.
			The flow assigned to $C_m$ also follows the same rules described in Proposition ~\ref{prop:subflow_even_wheel_41_edges} for the letter and the sign, but we reverse the rule for the presence of $\sim$: more precisely, if $m$ is odd we assign a flow without $\sim$ and if $m$ is even we assign a flow with $\sim$.  
			In this way, we guarantee that the number of components which reverse the orientation is even. Indeed, the number of odd fans plus the number of odd connectors is odd, but we recover the parity by the modification on $C_m$. Hence a $\sigma$-faithful flow is defined in $W_n$.
		\end{proof}

		\subsection{Even subgraphs with edges of capacity $(1,2)\cup(3,4)$}
		
		Until now we have given a complete characterization of graphs with circular flow number at least $5$ arising from wheels having an even subgraph of $(4,1)$-edges (and then also for $(4,0)\cup (0,1)$-edges). Note that several methods in~\cite{EMT16} concerned edges with capacity of measure $2$. Hence it is natural to ask whether a similar result holds for even subgraphs with the edges of measure 2 we are left with, i.e.\ $(1,2)\cup(3,4)$-edges. In this section, we  prove that the behavior is slightly different in this case.
		
		Select, for the entire section, $x,y\in (1,2)\cup(3,4)$ such that $x=1+2\delta$ and $y=3.5+2\delta$ with $\delta \in (0,0.25)$. Note that the difference between $x$ and $y$ is exactly $2.5$ for every choice of $\delta$.
		
		First of all, we completely solve the case where the even subgraph $J$ is induced by the edges of the external cycle of the wheel $W_n$.
		
		\begin{theorem}\label{teo:wheel_12_ext_cycle}
			Let $G$ be a graph in $(W_n,J_{(1,2)\cup(3,4)})$ where $J$ is the even subgraph induced by the edges of the external cycle of the wheel $W_n$. Then
			\[ \Phi_c (G)\ge 5 \text{ if and only if } n \text{ is odd. } \]
		\end{theorem}
		
		\begin{proof}
			Suppose $n$ is odd. Then, the assertion follows as a direct application of Corollary~\ref{Cor:P_odd_cycle_1}.
			Suppose $n$ is even. Now, take an orientation of $W_n$ such that the external cycle of $W_n$ is clockwise oriented, and  assign alternately flow values $x$ and $y$ to the edges of the external cycle. Finally, assign flow value $2.5$ to all edges incident to the central vertex (since the flow value is $2.5$ modulo $5$ the orientation of such edges does not really matter). The defined flow in $W_n$ is a $\sigma$-faithful flow, then $\Phi_c(G)<5$ for every $G\in (W_n,J_{(1,2)\cup(3,4)})$ by Proposition~\ref{pro:Hsigma}.
		\end{proof}
		
		Now, we consider the more general case in which $J$ is an arbitrary even subgraph of $W_n$. We use the terminology introduced in the previous section to describe $J$. In order to characterise all even subgraphs $J$ such that $(W_n, J_{(1,2)\cup(3,4)})$ contains graphs with circular flow number at least 5, we need to define some particular flows on $l$-fans and $m$-connectors.
		
		\underline{Flows $f_x$ and $f_y$ in an $l$-fan:} consider an $l$-fan $F_l$, we assign a clockwise orientation to its external cycle of $(1,2)\cup(3,4)$-edges and we assign alternately the values $x$ and $y$ to each $(1,2)\cup(3,4)$-edge. 
		Note that, if $l$ is even, the first edge and the last edge of $F_l$ receive the same flow value, otherwise, if $l$ is odd, they receive distinct values. 
		Finally assign to all further edges, both internal and lateral, the flow value $2.5$. Once again, note that the orientation of edges with flow value $2.5$ is not relevant, so we can choose lateral edges to be oriented in a clockwise direction. 
		It is important to note that if we add the flow value $2.5$ on the external cycle of the $l$-fan in clockwise direction, we obtain a new flow having edges with flow values $x$ and $y$ exchanged.
		We denote by $f_x$ and $f_y$ these two flows on an $l$-fan. More precisely, $f_x$ and $f_y$ are the flow defined as above and having flow value $x$ and $y$, respectively, on the first edge of $F_l$.
		
		\underline{Flow $g_x$ and $g_y$ in an $m$-connector, $m\neq1$:} Consider an $m$-connector, with $m\ge2$, we proceed exactly as we did above for an $l$-fan: we assign a clockwise orientation to the external cycle and, since $(1,2)\cup(3,4) \subset (1,4)$, we can assign flow values exactly as for an $l$-fan.
		Again, we use the notation $g_x$ and $g_y$ to denote the flows having values $x$ and $y$, respectively, on the unique edge of the external cycle of the $m$-connector directed away from $v_c$. 
		Finally, if  $m=0$, we simply give to the unique edge a clockwise orientation with respect to the external cycle of $W_n$ and we assign flow value $2.5$: we will denote such a flow by $g_x$.
		
		Previous flows are defined in $l$-fans, for any possible $l$, and $m$-connectors, for $m\neq1$.
		In order to deal with $1$-connectors we are going to define two methods to obtain a flow in them. 
		Both methods partially affect the flow values on some edges of the fans adjacent to the $1$-connector, but, with a suitable choice of the parameters, the resulting flows are still $\sigma$-faithful flows.
		
		Take $x$ and $y$ as described before, and set $\delta'=0.5+\delta$.
		 
\begin{figure}[htb!]
			\centering
			\includegraphics[width=8cm]{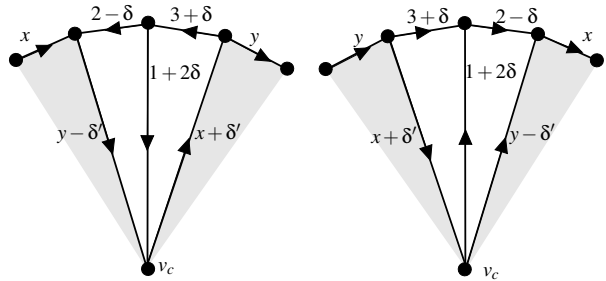}
			\caption{Method A to assign a flow to a $1$-connector.}\label{Fig:Meth_A}
		\end{figure}

		\underline{Method A:} Consider a $1$-connector $C$ and let $F$ and $F'$ be the two fans which precede and follow $C$, respectively. Assume that flows $f_x$ or $f_y$ are assigned on $F$ and $F'$ in such a way that the last edge of $F$ and the first edge of $F'$ have different flow values. That is: one of them has flow value $x$ and the other has flow value $y$ -- the two possible cases are presented in Figure~\ref{Fig:Meth_A}. Then, according to the flows on $F$ and $F'$, we can modify the flow value of the last edge of $F$ and the first edge of $F'$ as in Figure~\ref{Fig:Meth_A}. Moreover, the same figure shows a way to assign a suitable orientation and flow value to all edges of $C$. The result is a new zero-sum flow for all vertices of the external cycle of the wheel and each flow value belongs to the interval assigned by the capacity function $\sigma$. 
		
\begin{figure}[htb!]
			\centering
			\includegraphics[width=9cm]{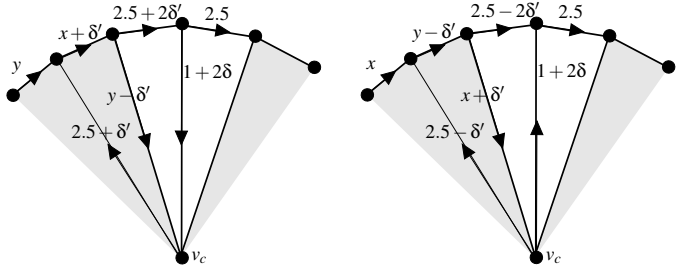}
			\caption{Method B to assign a flow to a $1$-connector.}\label{Fig:Meth_B}
		\end{figure}

		\underline{Method B:} Consider a $1$-connector $C$ and let $F$ be an $l$-fan, with $l>2$, adjacent to $C$. Assume that a flow $f_x$ or $f_y$ is assigned on $F$: the two possible cases are presented in Figure~\ref{Fig:Meth_B}. Then we can modify the flow value of three edges of $F$ and we can assign an orientation and a flow value to all the edges of $C$ as shown in Figure~\ref{Fig:Meth_B}. The result is a zero-sum flow for all vertices of the external cycle of the wheel and each flow value belongs to the interval assigned by the capacity function $\sigma$.
			
		\begin{theorem}\label{teo:wheel_12_not_ext_cycle}
			If $G\in(W_n,J_{(1,2)\cup(3,4)})$, where $J$ is not the external cycle and $n\geq4$, then $\Phi_c(G)<5$.		
			If $n=3$, $\Phi_c(G)<5$ if and only if $J$ is not a $3$-cycle.
		\end{theorem}
		
\begin{figure}[h]
			\centering
			\includegraphics[width=5cm]{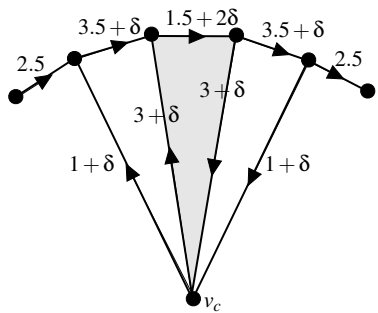}
			\caption{The flow assigned to a $2$-fan and its adjacent $1$-connectors in the proof of Theorem~\ref{teo:wheel_12_not_ext_cycle}.}\label{Fig:sp2fan}
		\end{figure}

		\begin{proof}
			
			If $n=3$ and $J$ is a $3$-cycle of $W_3$, then the unique vertex of $W_3$ not in $J$ can always be viewed as the central vertex and $J$ as the external cycle of $W_3$. In this case, the result follows by Theorem~\ref{teo:wheel_12_ext_cycle}. Otherwise, if $J$ is a $4$-cycle of $W_3$, we can assign the flow $g_x$ to the unique $0$-connector and the flow $f_x$ to the unique $3$-fan, thus obtaining the required flow.
			
			Now, consider $n\geq4$ and $J$ an even subgraph of $W_n$ distinct from the external cycle of $W_n$. We will prove that for any possible $J$ there exists a $\sigma$-faithful flow in $W_n$. Then, the assertion will follow by Proposition~\ref{prop:circ_F_N_sub_Flow} and~\ref{pro:Hsigma}.  
			In each step of the proof we construct the $\sigma$-faithful flow by first assigning flows in $l$-fans and $m$-connectors with $m\neq1$, and then by applying Method A and Method B to assign a flow in $1$-connectors as well.
			
			\emph{CLAIM 1:} If $J$ has an $m$-connector $C$ with $m\neq1$, then $\Phi_c(G)<5$.
			
			\emph{Proof of Claim 1:} Let $F$ be the fan following $C$ in $J$. Orient every fan in $J$ in clockwise direction. Starting from the first edge of $F$, we follow the assigned orientation on the edges of $J$ and we alternately assign flow value $x$ and $y$ to its edges. Assign an arbitrary orientation and flow value $2.5$ to all internal edges of every fan of $J$. In this way, we have defined in each fan a flow which is either $f_x$ or $f_y$, according to the flow value of the first edge of the fan.
			Now we define the flow in the connectors. We assign the flow $g_x$ to every $m$-connector with $m\neq1$ (connector $C$ included). Finally, we use Method A to assign the flow to every $1$-connector. Note that, in this case, we can apply Method A to every $1$-connector because $C$ is the unique connector such that flow values of the last edge of the fan preceeding it and the flow value of the first edge of the fan following it could be equal. 
			Hence, we have constructed a flow $f$ in $W_n$ with the required properties, then the claim follows by Proposition~\ref{pro:Hsigma}.\qed
			
			Hence, from now on, we can assume that all $m$-connectors of $J$ are $1$-connectors.
			
			\emph{CLAIM 2:} If $J$ has an $l$-fan $F$, with $l>2$, then $\Phi_c(G)<5$.
			
			\emph{Proof of Claim 2:} Let $C$ be the $1$-connector of $J$ which follows $F$ in $J$, and $F'$ the fan following $C$. Orient every fan in $J$ in clockwise direction. Starting from the first edge of $F'$, we follow the assigned orientation on the edges of $J$ and we alternately assign flow value $x$ and $y$ to its edges. Assign an arbitrary orientation and flow value $2.5$ to all internal edges of every fan of $J$. In this way, we have again defined a flow in each fan which is either $f_x$ or $f_y$, according to the flow value of the first edge of the fan.
			Now we define the flow on $1$-connectors. We can use Method A to assign the flow to every $1$-connector except, possibly, to $C$. Indeed, $C$ is the unique $1$-connector for which the flow value of the last edge of the previous fan could be equal to the flow value of the first edge of the next fan. Anyway, since $F$ is an $l$-fan with $l>2$, we can apply Method B to obtain a flow in $C$.
			Hence, we have constructed a $\sigma$-faithful flow in $W_n$, then the claim follows by Proposition~\ref{pro:Hsigma}.\qed
			
			Hence, from now on, we can assume that all connectors are $1$-connectors and all fans are $2$-fans.
			In order to complete the proof, we have to distinguish between two cases according to the parity of the number of $2$-fans in $J$.
			Assume that $J$ has an even number of $2$-fans (and then also an even number of $1$-connectors). Orient every fan in $J$ in clockwise direction. Starting from the first edge of an arbitrary $2$-fan, we follow the assigned orientation on the edges of $J$ and we alternately assign flow value $x$ and $y$ to its edges. Assign an arbitrary orientation and flow value $2.5$ to all internal edges of every fan of $J$. In this way, we have again defined a flow in each fan which is either $f_x$ or $f_y$, according to the flow value of the first edge of the fan.
			Now, we can use Method A to assign the flow to every $1$-connector since the flow value of the last edge of a fan is always different from the flow value of the first edge of the next fan in the sequence. Hence, we have constructed a $\sigma$-faithful flow in $W_n$ also in this case.
			Now, assume that $J$ has an odd number of $2$-fans (and then also an odd number of $1$-connectors). Select a $2$-fan $F$ of $J$. Assign a flow on $F$ and to the two connectors adjacent to $F$ as in Figure~\ref{Fig:sp2fan}. Assign alternately flow $f_x$ and $f_y$ to all other $2$-fans and use Method A to assign the flow to each $1$-connector between them. This defines a $\sigma$-faithful flow in $W_n$. 
		\end{proof}

\section{Even subgraphs with edges having a capacity set of measure different from 2}
\label{sect:even_subg}
	
	In Section~\ref{sect:main_results}, we have completely analysed some particular instances of the following general problem:
	
	\begin{problem}\label{main_problem}
		Given a wheel $W_n$ of length $n\ge3$ with a prescribed even subgraph $J$ and an element $A\in SI_5$, establish if $\Phi_c(G)\ge 5$ for $G \in (W_n,J_A)$.
	\end{problem}
	
	More precisely, we completely answered to all possible instances with $Me(A)=2$. Our goal in this section is to analyse all other cases when $Me(A)\neq 2$.
	
	Let us first recall the two methods mentioned in~\cite{EMT16} which make use of generalised edges of measure $0$ and $1$ capacity, respectively.
	
	\begin{itemize}
		\item[$M1$.] Let $G$ be a graph consisting of simple edges with a degree $3$ vertex $v$, then by replacing two of the edges adjacent to $v$ with $(2,3)$-edges we can generate a graph in $F_{\ge 5}$.
		\item[$M2$.] Insert an $\emptyset$-edge anywhere in a graph $G$. Clearly the resulting graph does not admit a sub-$5$-MCNZF.
	\end{itemize}
	
	We will refer to the first method as \emph{method M1} and to the second as \emph{method M2.}
	
	\subsection{Set $A$ of measure $0$}
			
			The unique set in $SI_5$ of measure 0 is obviously the empty-set.
			The unique method that involves $\emptyset$-edges is method $M2$. 
			Now we prove that this method indeed does not produce any new examples.
			
			\begin{theorem}
				Let $G_{uv}$ be an $\emptyset$-edge. Denote by $G'$ the (multi-)graph obtained from $G_{uv}$ by identifying $u$ and $v$, then $\Phi_c(G')\ge 5$.
			\end{theorem}
			
			\begin{proof}
				Suppose by contradiction that $\Phi_c(G')<5$. Then there is a sub-$5$-MCNZF $\phi$ in $G'$. Select an orientation in $G'$ such that all edges which arise from edges incident to $v$ in $G_{uv}$ are oriented towards $v$ and all edges which arise from edges incident to $u$ in $G_{uv}$ are oriented outward from $u$. Then $G_{uv}$ inherits an orientation and a flow from $G'$, that, with a slight abuse of terminology, we still call $\phi$, such that $\phi|_{E(G_{uv})} \subseteq (1,4)$ and \[\sum_{e\in E^+(u)} \phi(e) = \sum_{e\in E^-(v)} \phi(e) \mod{5}. \]
				If $x\in \R/5\Z$ is defined to be the common result of those summations, then $x\in CP_5(G_{uv}) = \emptyset$, a contradiction.
			\end{proof}
			
			This last result shows that whenever we generate a graph $G$ in $F_{\ge5}$ using method $M2$, then $G$ could also be generated by a suitable expansion of a smaller graph $H\in F_{\ge5}$, where $H$ is obtained by identifying the terminals of the $\emptyset$-edge that has been used to generate $G$.
			
		\subsection{Set $A$ of measure $1$}
			
			Similarly as for the case of measure $2$, we would like to present this case as an expansion of a suitable wheel. To this purpose let us call a \emph{wheel of length $2$}, denoted by $W_2$, a loopless (multi)graph with exactly $2$ vertices of degree $3$ and a vertex of degree $2$.
			
			Consider a pair $(H,\sigma)$ that presents the configuration described in method $M1$. Call $v$ the vertex of degree $3$ and let $N_H(v)=\{w_1,w_2,w_3\}$ with both $w_1v$, $w_2v$ $(2,3)$-edges and $w_3v$ a simple edge. Identify all vertices in $V(H)-\{v\}$ to a unique vertex $w$, thus obtaining a multigraph with two vertices, $v$ and $w$, and three parallel edges between them, two of them are $(2,3)$-edges and one of them, say $e$, is a simple edge.
			
			Now, if we subdivide the unique simple edge $e$ with a new vertex, then we obtain a wheel of length $2$ with the external $2$-cycle consisting of $(2,3)$-edges. The subdivision operation does not alter the circular flow number of the graph because it generates a degree $2$ vertex that we can get rid of at any time by using the smoothing operation.
			
			Also note that every wheel having an even subgraph of $(2,3)$-edges presents the configuration described in $M1$, and so it can be reduced by contraction to a wheel of length $2$. Moreover, the presence of the configuration described in method $M1$ assures that a graph produced in this way has circular flow number $5$ or more.
			
			We would like to stress that this is the unique case in which a wheel of length $2$ does produce examples of graphs with circular flow number at least $5$ by using methods described in this paper.
			
	\subsection{Set $A$ of measure $3$}
			
			\begin{lemma}\label{lem:sub_flow_in_g-graphs}
				Consider the families $G^\sigma$ and $G^\rho$, such that $\rho(e)\subseteq \sigma(e)$ and $\Me(\rho(e))=\Me(\sigma(e))$ for every $e\in E$. Let $H_1\in G^\sigma$ and $H_2\in G^\rho$. Then a sub-$5$-MCNZF exists in $H_1$ if and only if it exists in $H_2$.
			\end{lemma}
			
			\begin{proof}
				The thesis follows from the fact that a sub-$5$-MCNZF can be taken with no integer values, possibly after adding a small quantity $\epsilon >0 $ to suitable directed cycles. 
			\end{proof}
			
			It follows by Proposition~\ref{pro:A'_subset_A} and Lemma~\ref{lem:sub_flow_in_g-graphs} that we can restrict our analysis of sets of measure $3$ to the sets $(1,4)$ and $(4,1)\cup (2,3)$.
			
			Note that an attempt to characterise graphs in $F_{\ge5}$ with only simple edges and $(1,4)$-edges is equivalent to ask for a  direct characterization of $F_{\ge5}$.
			Obviously, there is no wheel with all edges of capacity $(1,4)$ having circular flow number at least $5$, since it is a planar graph.
			
			Hence, we can focus on $(4,1)\cup (2,3)$-edges.
			
			\begin{theorem}\label{teo:measure_3_wheel}
				For every $n\in\N$, if $G\in(W_n,J_{(4,1)\cup (2,3)})$, then $\Phi_c(G)<5$.
			\end{theorem}
			
			\begin{proof}
				Since $(4,1) \subseteq (4,1)\cup (2,3)$, every $\sigma$-faithful flow presented in the case of $(4,1)$-edges is a $\sigma$-faithful flow also in this case. Hence, by Theorem~\ref{teo:odd_wheels_41_edges}, $\Phi_c(G)<5$ when either $n$ is even or $n$ is odd and there is no $k$-connector with $k>1$.
				
				So let us assume that $n$ is odd and $G$ only contains $0$ and $1$-connectors.
				
				Let $F$ be an $l$-fan and apply the procedure described in Proposition~\ref{prop:subflow_even_wheel_41_edges} starting from $F$, except for the first edge of $F$.	
				Recall that $z+x$ is in $(0,1)$. Then, there are suitable values of $x\in (1,2)$ and $z\in(4,0)$ such that $z+x+x\in(2,3)\subseteq (4,1)\cup (2,3)$.
				In order to complete the description of the flow, we orient $e$ away from $v_c$ and  assign it the flow value $z+x+x$. One can easily check that the defined flow is a zero-sum flow at every vertex and so it is a $\sigma$-faithful flow.
				
				The last case to be discussed is the one where the length of the wheel is odd and the chosen even subgraph is the external cycle $x_1\dots x_{2t+1}$, for $t\geq 2$. 
				
				Here, we define a $\sigma$-faithful flow $f$ as follows: we orient  the external cycle of the wheel in clockwise direction and we orient  every edge $v_cx_i$ with $i$ odd except for $i=1$ towards the central vertex $v_c$, and away from $v_c$ otherwise.
				
			We define flow values as follows: 
			\begin{itemize}
				\item $f(x_ix_{i+1}):=z$ for $i$ odd $\in \{3,\dots, 2t-1\};$
				\item $f(x_ix_{i+1}):=z+x$ for $i$ even $\in\{4,\dots, 2t\};$
				\item $f(x_1x_{2t+1}):= z;$
				\item $f(x_1x_2):= z+x;$
				\item $f(x_2x_3):=z+x+x;$
				\item $f(v_cx_3):=x+x$ and $f(v_cx_i):=x$ for $i \neq 3$.
			\end{itemize}
			
			Therefore there is no graph with circular flow number at least $5$ belonging to the family $(W_n,J_{(4,1)\cup (2,3)})$.
\end{proof}
		
			\subsection{Set $A$ of measure $4$ or $5$}
		
			If the set $A$ has measure either $4$ or $5$, it will turn out that no instance produces a graph in $F_{\ge5}$.
			
			\begin{theorem}
			For every $n\in\N$ and every $A \in SI_5$ with $Me(A)\geq4$, if $G\in(W_n,J_A)$, then $\Phi_c(G)<5$.
			\end{theorem}
			\begin{proof}
			If $Me(A)=5$, we can simply observe that each open integer set of measure $5$ contains $(4,0)\cup (0,1)\cup(2,3)$, for which Theorem~\ref{teo:measure_3_wheel} holds.
			
			If $Me(A)=4$ then Lemma~\ref{lem:sub_flow_in_g-graphs} says that we can only consider the case $A=(3,2)$. Since $(1,2)\cup(3,4)$ is a subset of $(3,2)$, from Theorem~\ref{teo:wheel_12_ext_cycle} and~\ref{teo:wheel_12_not_ext_cycle} we deduce that the only case that could produce a graph in $F_{\ge5}$ is when the even subgraph $J$ is the external cycle and the wheel has odd length. But now we will show that also in this case we do not obtain any graph in $F_{\ge5}$.
			
			Consider an odd wheel $W_{2t+1}$ and let $J$ be the even subgraph induced by the edges of the external cycle $x_1x_2\dots x_{2t+1}$. Assume $\sigma(e)=(3,2)$ for every edge $e$ of $J$. 
			Let $x,z\in (1,2)$ such that $y:=x+z\in (3,4)$. There exists an $\alpha>1$, sufficiently close to $1$, such that \begin{itemize}
				\item $\alpha + z \in (1,4)$.
				\item $y-\alpha - z = x - \alpha \in (0,1)$. 
			\end{itemize} 
			
			Hence we can define a $\sigma$-faithful flow $f$ in $W_{2t+1}$.
			We orient the external cycle of the wheel in clockwise direction and we orient  every edge $v_cx_i$ with $i$ odd except for $i=1$ towards the central vertex $v_c$, and away from $v_c$ otherwise.
			We define flow values as follows: 
			\begin{itemize}
				\item $f(x_ix_{i+1}):=x$ for $i$ odd $\in \{1,\dots, 2t-1\};$
				\item $f(x_ix_{i+1}):= y $ for $i$ even $\in\{2,\dots, 2t\};$
				\item $f(x_{2t+1}x_1):= y-\alpha-z;$
				\item $f(v_cx_i):=z$ for $i \in\{2,\dots,2t\};$
				\item $f(v_cx_{2t+1}):= \alpha - z$ and $f(v_cx_1):= \alpha$.
			\end{itemize}
			
			Therefore the cases $k=4$ and $k=5$ do not produce any new examples of graphs with circular flow number $5$ or more.

			\end{proof}

		\section{Final remarks}
		\label{sect:final_remarks}
		
		In this paper we have considered all possible instances of Problem~\ref{main_problem} and we have proved that several known constructions of graphs with circular flow number at least 5 can be described as particular instances of this problem.  
		
		All our results can be summarised in the following theorem:
		
		\begin{theorem}\label{theo:summary}
		If $G \in (W_n,J_I)$ where $J$ is a (non-empty) even subgraph of $W_n$ and $I \in SI_5$, then $\Phi_c(G)\geq5$ if and only if one the following holds
		\begin{itemize}                                                                                                                   \item $I=\emptyset$;                                                                                                    \item $I=(2,3)$;
		\item $I \subseteq (4,1)$, $n$ odd and $J$ has no $k$-connector for $k>1$;
		\item $I \subseteq (1,2) \cup (3,4)$ and either $n>3$ odd and $J$ is the external cycle of $W_n$ or $n=3$ and $J$ is a $3$-cycle.        
		\end{itemize}

		\end{theorem}

		A more general formulation of the previous problem could be considered, where the capacity function $\sigma$ is not constant in $J$.
		
		\begin{problem}\label{pro:general}
		Given a wheel $W_n$ with $n+1$ vertices and $J$ a non-empty even subgraph of $W_n$, establish for every integer $n$, every even subgraph $J$ and every capacity function $\sigma$ if a graph $G \in W_n^\sigma$ has circular flow number at least 5, where $\sigma(e)=(1,4)$ for all $e \notin E(J)$.  
		\end{problem}
		
		In particular, our results say that all methods from~\cite{EMT16} can be described as an instance of Problem~\ref{main_problem}, i.e.\ considering $\sigma$ constant on all edges of $J$, except the one arising from Lemma~4.6 in~\cite{EMT16} where we need to consider the more general Problem~\ref{pro:general}; indeed, in this last case, we must consider a capacity function $\sigma$ which could be not constant on $J$.
		
		Last observation suggests that new methods can be obtained by looking at Problem~\ref{pro:general} in its general formulation, and we leave this as a possible research problem.
		
		It is interesting to note that nearly all small snarks with circular flow number $5$ can be obtained by using methods described in this paper. 
		M\'{a}\v{c}ajov\'{a} and Raspaud determined all snarks with circular flow number $5$ up to 30 vertices in~\cite{MaRa06}. We designed an algorithm for computing the circular  flow number of a cubic graph (the details of this algorithm will be described in~\cite{GMM2}). By applying this algorithm to the complete list of all snarks up to 36 vertices from~\cite{BGHM}, we were able to determine all snarks with circular flow number 5 up to that order. The counts of these snarks can be found in Table~\ref{Table:snarks}.

	\begin{table}[htb!]
		\centering 
		 \renewcommand{\arraystretch}{1.1} 
		\begin{tabular}{|c || r | c |} 
			\hline 
			Order & Number of snarks & Number of snarks in $F_{\ge5}$  \\ 
			\hline 
			10 & 1 & 1 (1) \\ 	
			12 & 0 & 0 (0) \\ 	
			14 & 0 & 0 (0) \\ 	
			16 & 0 & 0 (0) \\ 	
			18 & 2 & 0 (0) \\ 	
			20 & 6 & 0 (0) \\ 	
			22 & 20 & 0 (0) \\ 	
			24 & 38 & 0 (0) \\ 	
			26 & 280 & 0 (0) \\ 		
			28 & 2900 & 1 (1) \\
			30 & 28~399 & 2 (2) \\
			32 & 293~059 & 9 (9) \\
			34 & 3~833~587 & 25 (25) \\
			36 & 60~167~732 & 98 ({\bf 96}) \\
			\hline 
		\end{tabular}
		\caption{Counts of all snarks with circular flow number $5$ up to 36 vertices. We indicate the number of snarks with circular flow number $5$ which can be obtained by using methods of the present paper in parentheses.}\label{Table:snarks}
	\end{table}		
			
		We also verified which snarks with circular flow number $5$ fit our description which led to the following observation (see also Table \ref{Table:snarks}).
		
\begin{observation}\label{obs:cfn_fit}
All snarks with circular flow number $5$ with order at most $34$ and  $96$ out of the $98$ snarks with circular flow number $5$ with order $36$ can be obtained by our unified construction method from Section~\ref{sect:constructions}.
\end{observation}		


All graphs from Table~\ref{Table:snarks} can be downloaded from the \textit{House of Graphs}~\cite{hog} at \url{http://hog.grinvin.org/Snarks}. The snarks with circular flow number $5$ can also be inspected at the database of interesting graphs from the \textit{House of Graphs} by searching for the keywords ``snark with circular flow number 5''. The two snarks with circular flow number $5$ on 36 vertices from Observation~\ref{obs:cfn_fit} which cannot be obtained through our unified method can be found by searching for ``snark with circular flow number 5 which cannot be obtained''. These two snarks will also be discussed in~\cite{GMM2}.
		
		All snarks obtained by using our methods as well as the two additional snarks with circular flow number 5 on 36 vertices are cyclically $4$-edge-connected. Therefore the following open problem remains.
		
\begin{problem}
Is the Petersen graph the only cyclically $5$-edge-connected snark with circular flow number $5$?
\end{problem}		

		\section*{Acknowledgements}
		All computations for this work were carried out using the Stevin Supercomputer Infrastructure at Ghent University.

\section*{Appendix}

For the reader's convenience we here present an example of a construction of a snark from~\cite{EMT16} which makes use of all considerations introduced in Section~\ref{sect:generalised_edges}.
	
\begin{figure}[htb!]
			\begin{center}
				\begin{minipage}[c]{.40\textwidth}
					\centering
					\includegraphics[scale=0.7]{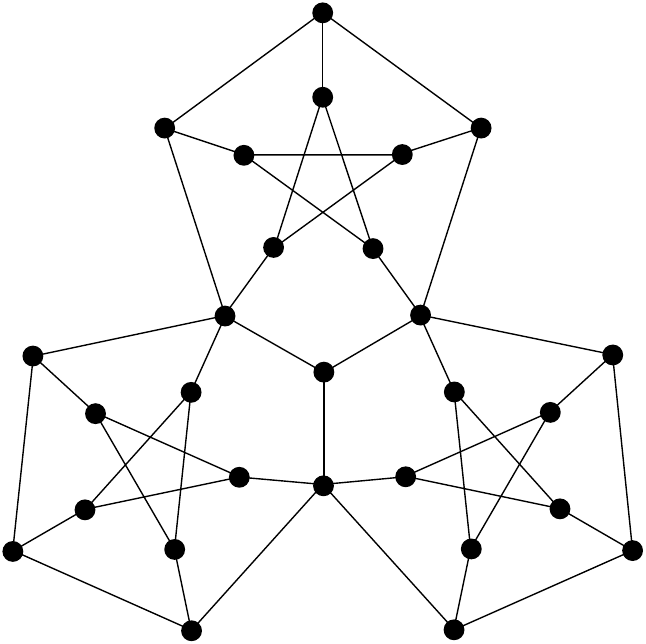}
				\end{minipage}%
				\hspace{4mm}%
				\begin{minipage}[c]{.40\textwidth}
					\centering
					\includegraphics[scale=0.7]{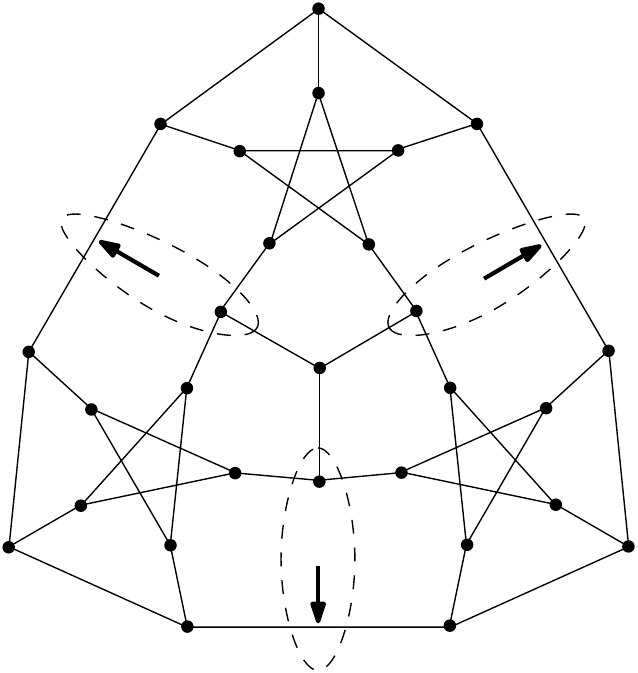}
				\end{minipage}
			\end{center}
			\caption{Expansion of the degree $5$ vertices of a wheel $W_3$ having an external cycle of $(4,1)$-edges.}\label{Fig:Espansione}
		\end{figure}
		
Denote by $C_3$ a $3$-cycle of $K_4$. Consider the pair $(K_4, \sigma)$ where $\sigma(e)=(4,1)$ for all $e \in C_3$ and $\sigma(e)=(1,4)$ otherwise. Thanks to Corollary~\ref{Cor:P_odd_cycle_1} we can say that all graphs in $K_4^{\sigma}$ have circular flow number at least 5.	In particular, if we replace every edge of $C_3$ with the generalised edge $\ca P^*_{10}(u,v)$, that is the Petersen graph minus an edge, of capacity $(4,1)$. We obtain a graph $G \in K_4^{\sigma}$ with circular flow number at least $5$. The graph $G$ is not a snark as it is not cubic yet. 
Now, we expand each vertex of degree $5$ of $G$ to the graph with two isolated vertices, and we connect these new vertices to the rest of the graph as in Figure~\ref{Fig:Espansione}. This operation produces three vertices of degree $2$ and we remove them by smoothing. None of these operations reduces the circular flow number. Note that the final graph has girth at least $5$ and is cyclically $4$-edge-connected, so it is a snark.			
Note also that such a snark is the smallest one larger than the Petersen graph and having circular flow number 5 (see~\cite{MaRa06} and Table~\ref{Table:snarks}).
		

\begin{thebibliography}{99}

 
\bibitem{AKLM} M. Abreu, T. Kaiser,  D. Labbate, and G. Mazzuoccolo, Tree--like snarks, {\em Electron. J. Combin.}, 23 (2016), no. 3, Paper 3.54.

\bibitem{hog}
G.~Brinkmann, K.~Coolsaet, J.~Goedgebeur, and H.~M{\'e}lot, House of Graphs: a database of interesting graphs, {\em Disc. Appl. Math.}, 161 (2013), 311--314. Available at \url{http://hog.grinvin.org/}.
 
\bibitem{BGHM}
 G.~Brinkmann, J.~Goedgebeur, J.~H\"agglund, and K.~Markstr\"om, Generation and properties of snarks, {\em J. Combin. Theory Ser. B.}, 103 (2013), 468--488.
  
 \bibitem{EMT16} L. Esperet, G. Mazzuoccolo, and M. Tarsi, The structure of graphs with Circular flow number 5 or more, and the complexity of their recognition problem, {\em J. of Combinatorics}, 7 (2016), 453--479.

\bibitem{GTZ} L.A. Goddyn, M. Tarsi, and C.Q. Zhang, On $(k,d)$-colorings and fractional nowhere-zero flows, \emph{J. Graph Theory}, 28 (1998), 155--161. 


\bibitem{GMM2} J.~Goedgebeur,  D. Mattiolo, and G. Mazzuoccolo, On the circular flow number of snarks, \emph{in preparation}.
 
\bibitem{LuSk} R. Lukot'ka and M. {\v{S}}koviera, Snarks with given real flow
numbers, \emph{J. Graph Theory}, 68 (2011), 189--201.
 
\bibitem{MaRa06} E. M\'{a}\v{c}ajov\'{a} and A. Raspaud, On the strong circular 5-flow conjecture, {\em J. Graph Theory}, 52 (2006), 307--316.

\bibitem{Mo} B. Mohar, Problem of the Month, \url{ http://www.fmf.uni-lj.si/~mohar/Problems/P03034CircularFlowConjecture.html}, March and April 2003.

\bibitem{PaZh} Z. Pan and X. Zhu, Construction of Graphs with Given Circular Flow
Numbers,  \emph{J. Graph Theory}, 43 (2003), 304--318.

\bibitem{SS} M. Schubert and E. Steffen, The set of circular flow numbers of regular graphs,  {\em J. Graph Theory}, 76 (2014), 297--308.

\bibitem{S} E. Steffen, Circular flow numbers of regular multigraphs, {\em J. Graph Theory}, 36 (2001), 24--34.

\bibitem{Tu} W.T. Tutte, A contribution to the theory of chromatic
  polynomials, \emph{Canad. J. Math.}, 6 (1954), 80--91.

\bibitem{Z} C.Q. Zhang, Integer Flows and Cycle Covers of Graphs, \emph{New
York: Marcel Dekker} (1997).  
  
\end{thebibliography}
\end{document}